\documentclass[12pt,a4paper]{article}
\usepackage[]{fontenc}
\usepackage[latin1]{inputenc}
\usepackage[francais]{babel}
\usepackage{a4,amsmath,amssymb,amsfonts,amsthm}
\usepackage{amsrefs}
\usepackage[all]{xy}

\begin{document}
\bibliographystyle{plain}
 
% Diff\'erentes fontes

\def\MF#1{\mathbb#1}     % fontes tableau noir (ensembles R, C, Z, ...)
\def\BF#1{\bold#1}       % symboles math\'ematiques accentues
\def\CF#1{\mathcal#1}    % lettres caligraphiees (mode math ou texte)
\def\EF#1{\scr#1}        % autres fontes caligraphiees (mode math ou texte)

% lettres spéciales

 \def\mQ{\MF{Q}}
 \def\mR{\MF{R}}
 \def\mZ{\MF{Z}}
 \def\mC{\MF{C}}
 \def\mN{\MF{N}}
 \def\mX{\MF{X}}
\def\mP{\MF{P}}
\def\mF{\MF{F}}
 \def\CI{{\cal I}}
 \def\CH{{\cal H}}
 \def\CO{{\cal O}}
 \def\CA{{\cal A}}
 \def\CB{{\cal B}}
 \def\CV{{\cal V}}
\def\CX{{\cal X}}
 \def\CC{{\cal C}}
\def\CN{{\cal N}}
\def\CL{{\cal L}}
\def\CM{{\cal M}}
\def\CF{{\cal F}}
\def\CD{{\cal D}}

% roman dans $$

\def\Spec{{\rm Spec}}
\def\rg{{\rm rg}}
\def\rHom{{\rm Hom}}
\def\Aut{{\rm Aut}}
 \def\Tr{{\rm Tr}}
 \def\Exp{{\rm Exp}}
 \def\Gal{{\rm Gal}}
 \def\End{{\rm End}}
 \def\det{{{\rm det}}}
 \def\Td{{\rm Td}}
 \def\ch{{\rm ch}}
 \def\che{{\rm ch}_{\rm eq}}
  \def\Spec{{\rm Spec}}
\def\Id{{\rm Id}}
\def\Zar{{\rm Zar}}
\def\Supp{{\rm Supp}}
\def\eq{{\rm eq}}
\def\Ann{{\rm Ann}}
\def\rad{{\rm rad}}
\def\Jac{{\rm Jac}}
\def\char{{\rm char}}
\def\lng{{\rm long}}
\def\rN{{\rm N}}
 \def\rh{{\rm h}}
 \def\rM{{\rm M}}
 \def\rd{{\rm d}}
 \def\rc{{\rm c}}
 \def\rH{{\rm H}}
 \def\rB{{\rm B}}
 \def\Bu{{\rm Bu}}
 \def\rK{{\rm K}}
 \def\rR{{\rm R}}
 \def\rL{{\rm L}}
 \def\re{{\rm e}}
 \def\rv{{\rm v}}
 \def\Tor{{\rm Tor}}
 \def\rH{{\rm H}}
 \def\Pro{{\rm Pro}}
\def\Ind{{\rm Ind}}
\def\Ens{{\rm Ens}}
\def\Pic{{\rm Pic}}
\def\Max{{\rm Max}}
\def\Sup{{\rm Sup}}
\def\mod{{\rm mod}}
\def\ordre{{\rm ordre}}
\def\Exc{{\rm Exc}}
\def\GL{{\rm GL}}

% général 

 \def\beginProof{\par{\bf Preuve. }}
 \def\endProof{${\qed}$\par\smallskip}
 \def\mtr#1{\overline{#1}}
 \def\trelp{{\rm trelp}}
\def\trel{{\rm trel}}
 \def\unr{{\rm unr}}
 \def\refeq#1{(\ref{#1})}
 \def\blb{{\big(}}
 \def\brb{{\big)}}
\def\mc{{{\mathfrak c}}}
\def\mcpr{{{\mathfrak c}'}}
\def\mcprpr{{{\mathfrak c}''}}
\def\ul#1{\overline{#1}}
\def\ss{{\rm ss}}
\def\parf{{\rm parf}}
\def\P1{{{\bf P}^1}}
\def\pr{\prime}
\def\prpr{\prime\prime}
\def\ss{\scriptstyle}
\def\OX{{ {\cal O}_X}}
\def\mpartial{{\mtr{\partial}}}
\def\inv{{\rm inv}}
\def\indlim{\underrightarrow{\lim}}
\def\prolim{\underleftarrow{\lim}}
\def\pprolim{'\prolim'}
\def\without{\backslash}
\def\pbdb{{\Pro_b\ D^-_c}}
\def\qc{{\rm qc}}
\def\Com{{\rm Com}}
\def\an{{\rm an}}
\def\gfield{{\rm\bf k}}
\def\sep{{\rm sep}}
\def\mod{{\rm mod}}
\def\mfp{{\mathfrak p}}
\def\mfl{{\mathfrak l}}
\def\mfn{{\mathfrak n}}
\def\adeg{\widehat{\rm deg}}
\def\ac1{\ari{\rc}_1}
\def\NT{{\rm NT}}
\def\wt#1{\widetilde{#1}}
\def\ari#1{\widehat{#1}}
\def\ul#1{\underline{#1}}
\def\c1{{\rm c}_1}
\def\og{$\scriptscriptstyle\langle\langle\ $}
\def\fg{$\ \scriptscriptstyle\rangle\rangle$}

 \newtheorem{theor}{Théorème}[section]
 \newtheorem{prop}[theor]{Proposition}
 \newtheorem{cor}[theor]{Corollaire}
 \newtheorem{lemme}[theor]{Lemme}
 \newtheorem{slemme}[theor]{sous-lemme}
 \newtheorem{defin}[theor]{Définition}
 \newtheorem{conj}[theor]{Conjecture}

 \parindent=0pt
 \parskip=5pt
 
 \author{
 Damian RÖSSLER\protect\footnote{Institut de Math\'ematiques de Jussieu,
 Universit\'e Paris 7 Denis Diderot, C.N.R.S.,
 Case Postale 7012,
 2 place Jussieu,
 F-75251 Paris Cedex 05, France,
 E-mail : dcr@math.jussieu.fr}}
 \title{Sur la distance $p$-adique entre un point d'ordre fini et une courbe de 
 genre supérieure ou égale à deux
 \\ \vskip10pt
 \small\bf On the $p$-adic distance between a point of finite order and a curve 
of genus higher or equal to two}
 \date{}
 \maketitle
\begin{center}\small
{\bf Résumé}
\end{center}{\small 
 Soit $A$ une variété abélienne sur $\mC_p$ ($p$ premier) 
 et $V\hookrightarrow A$ une sous-variété fermée. La conjecture 
 de Tate-Voloch prédit que la distance $p$-adique entre un point de 
 torsion $T\not\in V(\mC_p)$ et la variété $V$ est bornée inférieurement par une 
 constante strictement positive. Cette conjecture 
 est démontrée par Hrushovski et Scanlon, lorsque $A$ a un modèle sur $\mtr{\mQ}_p$. Dans le cas 
 où $V$ est une courbe plongée dans sa jacobienne et possède un modèle sur 
 un corps de nombres, nous donnons une formule 
 explicite pour cette constante, qui dépend d'invariants analytiques 
 et arakeloviens.}
\begin{center}\small
{\bf Résumé en anglais / English abstract}
\end{center}
{\small 
Let $A$ be an abelian variety over $\mC_p$ ($p$ a prime number) 
and $V\hookrightarrow A$ a closed subvariety. The conjecture 
of Tate-Voloch predicts that the $p$-adic distance from 
a torsion point $T\not\in V(\mC_p)$ to the variety $V$ is 
bounded below by a strictly positive constant. This conjecture 
is proven by Hrushovski and Scanlon, when $A$ has a model 
over $\mtr{\mQ}_p$. We give an explicit formula for this constant, 
in the case where $V$ is a curve embedded into its Jacobian and 
$V$ has a model over a number field. This explicit formula 
involves analytic and arakelovian invariants of the curve.}

\section{Introduction}

Soit $A$ une variété abélienne définie sur $\mC_p$, 
la complétion de la clôture algébrique du corps 
$\mQ_p$ des nombres $p$-adiques. Soit $V$ une sous-variété fermée 
de $A$. La {\it conjecture de Tate-Voloch} (cf. \cite[avant la Prop. 3]{Tate-Voloch-Linear}) prédit l'existence d'une constante 
$M_V>0$, avec la propriété suivante. Pour tout point de torsion 
$T\in\Tor(A(\mC_p))$, si $\rd_p(V,T)<M_V$ alors $T\in V(\mC_p)$. 
Ici $\rd_p(\cdot,\cdot)$ est la distance $p$-adique entre 
$T$ et $V$ (voir plus bas pour cette notion). Cette conjecture 
est démontrée par Hrushovski et Scanlon (\cite[Introduction]{Hrushovski-Manin} et 
\cite{Scanlon-The-conjecture}) lorsque $A$ a un modèle sur une extension finie de $\mQ_p$.

Le propos du présent texte est de donner une formule explicite pour la constante 
$M_V$, dans le cadre suivant: $V$ est définie sur un corps de nombres
 et l'immersion $V\hookrightarrow A$ est le plongement 
jacobien défini par un point rationnel. 
  Pour obtenir cette formule, nous combinons une technique de preuve 
de la conjecture de Tate-Voloch pour les courbes due à Buium (cf. plus bas) 
avec des résultats de Moret-Bailly en théorie d'Arakelov. 
Nous montrons que la constante $M_V$ peut être choisie comme 
une fonction d'invariants analytiques et arakeloviens (Th. \ref{mainth}). 
Dans le cas où la courbe a bonne réduction partout, ces invariants 
peuvent être approximés numériquement, si l'on dispose d'équations 
pour la courbe $C$ (cf. Prop. \ref{mainprop}). Par ailleurs, 
notre formule est uniforme en le nombre premier $p$ dans la mesure 
où les invariants sont indépendants du nombre premier $p$. 

L'idée de base de la démonstration est que le théorème du cube dans sa forme 
arakelovienne permet de comparer la distance $p$-adique 
d'un point de torsion au diviseur $\Theta$ avec la distance $p$-adique 
de l'origine au diviseur $\Theta$. Cette comparaison à elle seule 
donne une forme très faible de la conjecture de Tate-Voloch, 
dans la mesure où la distance $p$-adique 
y est bornée par une fonction du degré du point de torsion. 
Il reste alors à estimer ce degré; on l'estime par la procédure suivante. 
La méthode de Buium donne une estimée du degré résiduel 
du point de torsion; cette estimée implique une estimée 
pour l'ordre du point de torsion via les estimées de Hasse-Weil. 
Enfin l'estimée de l'ordre donne une estimée évidente pour le degré.

Venons en à la description précise de notre formule pour la 
constante $M_V$. 

Soit $C$ une courbe de genre $g\geqslant 2$, définie sur un corps 
de nombre $K_0$. Soit $P_0\in C(K_0)$ un point et 
$j=j_{P_0}:C\hookrightarrow\Jac(C)$ l'immersion jacobienne associée. 

Soit $K$ une extension finie de $K_0$ 
et soit $\mathfrak p$ un idéal premier dans $\CO_K$.  

On note $\CN(C_K)$ le modèle de Néron de $\Jac(C_K)$ sur $\CO_{K}$. 

Si $\mfl\in\Spec\ \CO_K$, on écrit $\Phi(\mfl)$ pour le groupe 
des composantes connexes de la fibre $\CN(C_K)_{\kappa(\mfl)}$ 
de $\CN(C_K)$ au-dessus de $\mfl$. 
On note $$\mfn_K={\rm ppcm}\{\#\Phi(\mfl)|\mfl\in\Spec\ \CO_K\}.$$

Soit $\wt{C}_\mfp$ la clôture de Zariski $\Zar(C_K)$ de $C_K$ sur dans $\CN(C_K)_\mfp:=
\CN(C_K)\times_{\CO_{K}}\CO_{K,\mfp}$. Ici $\CO_{K,\mfp}$ est la localisation de 
$\CO_K$ en $\mfp$. Soit $P\in \Jac(C)(K)=\CN(C_K)_\mfp(\CO_{K,\mfp})$. On définit alors
$$
\rv_\mfp(P,C_K):=\sup_{j\geqslant 0}\{j|{P}\ (\mod\ \mfp^j)\in\wt{C}_\mfp(\CO_K/\mfp^j)\}.
$$
On notera que $\rv_\mfp(P,C_K)$ vaut éventuellement $\infty$. \`A partir de là, 
on définit {\it la distance $\mfp$-adique} entre $P$ et $C$ par la formule 
$$
\rd_\mfp(P,C):=p^{-\rv_\mfp(P,C_K)/\re_{\mfp}}.
$$
Ici $\re_{\mfp}$ est le degré de ramification de $\mfp$ sur $K_0$.  
Par construction, la distance $\mfp$-adique est invariante 
par extension du corps de nombres $K$ et de la place $\mfp$. 

On définit la fonction de $m$
$$\Bu_m:=(m(2g-2)+6g)\cdot m^{2g}\cdot 3^g\cdot g!$$
et la fonction de $m$ et $n$
$$
\rL_{n,m}:=\big[n^{\Bu_m\cdot g}+(2^{2g}-2g-1)n^{\Bu_m(g-1)}+2gn^{\Bu_m(g-1/2)}\big]^{4g^2}.
$$
Enfin, on définit la fonction de $m$
$$
\rH_m:=L_{m^{[K_0:\mQ]},m}
$$
Si $Y\hookrightarrow\Jac(C)$ est un sous-variété fermée, on notera 
$$
\Exc(Y):=\Zar\big(Y(\mtr{K}_0)\cap\Tor(\Jac({C})(\mtr{K}_0))\big)
$$
Le sous-ensemble $\Exc(Y)\hookrightarrow Y_{\mtr{K}_0}$ est naturellement 
$\Gal(\mtr{K}_0|K_0)$-invariant et nous le considérerons donc comme un sous-schéma 
fermé réduit de $Y$. On rappelle que la conjecture de Manin-Mumford 
(théorème de Raynaud \cite{Raynaud-Sous-var}) prédit que les composantes irréductibles géométriques 
de $\Exc(Y)$ sont des translatées de sous-variétés abéliennes de 
$\Jac(C)_{\mtr{K}_0}$. 

Soit $S$ l'ensemble des nombres premiers que divise un idéal premier 
de mauvaise réduction de $C$ sur $K_0$. Soit 
$$
R:=2\cdot\CD_{K_0/\mQ}\cdot\mfn_{K_0}\cdot\prod_{l\in S}l
$$
où $\CD_{K_0/\mQ}$ est le discriminant de $K_0$ sur $\mQ$. 

Si $G$ est un groupe abélien et $n\in\mZ$, notons $\Tor^n(G)$ l'ensemble des éléments 
de $G$ qui sont d'ordre fini premier à $n$. 

\begin{theor}
On suppose que la courbe $C$ a réduction semi-stable sur $K_0$. Soit 
$L\subseteq\mtr{K}_0$ l'extension maximale 
de $K_0$ qui est non-ramifiée au-dessus des idéaux premiers 
de mauvaise réduction de $C$ sur $K_0$.

Soit
$$
D:=2\cdot[K_0:\mQ]\cdot\Big|\log\Theta_{\rm Max}(C_{\mtr{K}_0})+[K_0:\mQ]^{-1}\adeg(\mtr{\CO}(\Zar(W_{g-1}))_0)\Big|.
$$

Soit 
$p$ un nombre premier tel que $(p,R)=1$.  Soit $\ul{\mfp}\subseteq 
\CO_L$ un idéal premier divisant $p$. 

Alors l'inégalité
\begin{eqnarray*}
d_{\ul{\mfp}}(T,{C})&\geqslant&p^{\textstyle-\Big(1+ D\cdot \rH_p\Big)}
\end{eqnarray*}
est vérifiée pour tout $t\in\Tor^p(\Jac(C)(L))$ tel que 
$$
(g-1)T\not\in\Exc(W_{g-1})(\mtr{K}_0)\cup[-1]^*\Exc(W_{g-1})(\mtr{K}_0).
$$
\label{mainth}
\end{theor}

La proposition suivante relie l'invariant 
$[K_0:\mQ]^{-1}\adeg(\mtr{\CO}(\Zar(W_{g-1}))_0)$ apparaissant dans le Théorème 
\ref{mainth} à des invariants plus classiques.
\begin{prop}
Si $C$ a bonne réduction partout sur $K_0$, alors l'égalité
$$
[K_0:\mQ]^{-1}\adeg(\mtr{\CO}(\Zar(W_{g-1}))_0)={1\over 4}\NT(j(\omega_C))+
{1\over 2}\rh_{\rm Fal}(\Jac(C)_{\mtr{K}_0})+{g\over 4}\log(4\pi).
$$
est vérifiée.
\label{mainprop}
\end{prop}

{\it La fonction $\NT:A(\mtr{K}_0)\to\mR$ }
est la hauteur de Néron-Tate sur $\Jac(C)$  associée 
au diviseur $\Theta$ sur $\Jac(C)$. Le diviseur $\omega_C$ est 
le diviseur canonique de $C$. 

{\it L'invariant $\adeg(\mtr{\CO}(\Zar(W_{g-1}))_0)$} est un invariant 
arakelovien. Nous allons utiliser  les notations usuelles de la théorie 
d'Arakelov globale, comme décrites dans \cite[Par. 3.1]{Gillet-Soule-Arithmetic}. 
Nous allons travailler sur $\CO_{K_0}$, que nous considérons 
comme un anneau arithmétique, que nous munissons comme tel de tous  ses plongements dans $\mC$. 
Si $W$ est un schéma 
sur $\CO_{K_0}$, on écrira donc 
$$
W_\iota:=W\times_{\CO_{K_0},\iota}\mC
$$
pour tout plongement $\iota$ de $\CO_{K_0}$ dans $\mC$ et
$$
W_\mC:=\coprod_{\iota:\CO_{K_0}\hookrightarrow\mC}W_\iota
$$ 
(cf. \cite[Par. 3.2.1]{Gillet-Soule-Arithmetic}). 
On note $W_{g-1}:=C+C\dots+C\ ((g-1)\ {\rm fois})$. 
 On note  
$\Zar(W_{g-1})$ l'adhérence schématique de $W_{g-1}$ dans 
$\CN(C)$.  On note 
$\mtr{\CO}(\Zar(W_{g-1}))$ le fibré hermitien que l'on obtient 
en munissant ${\CO}(\Zar(W_{g-1}))$ de l'unique 
métrique $\langle\cdot,\cdot\rangle$ dont la forme de courbure est invariante par translation et telle que 
$$
\int_{\Jac(C)_\iota}\langle s,s\rangle\ \rd\mu_\iota=2^{-g/2},
$$
pour tout plongement $\iota$ de $K_0$ dans $\mC$. Ici $\rd\mu_\iota$ est la mesure de Haar sur $\Jac(C)_\iota$ de masse totale $1$ 
et $s$ est la  section canonique de ${\CO}(\Zar(W_{g-1}))$ 
(cf. \cite[Par. 3., (3.2.1)]{Moret-Bailly-Sur}). 
 On rappelle que l'opération 
$\adeg(\cdot)$ associe à chaque fibré en droites hermitien
sur $\Spec\ \CO_{K_0}$ un nombre réel; cf. 
\cite[Par. 2.1.2, (2.1.8)]{Bost-Gillet-Soule-Heights} pour la définition. 
On note enfin $\mtr{\CO}(\Zar(W_{g-1}))_0$ la restriction 
de $\mtr{\CO}(\Zar(W_{g-1}))$ à $\Spec\ \CO_{K_0}$ via la section nulle.
La définition $\adeg(\mtr{\CO}(\Zar(W_{g-1}))_0)$ est maintenant 
élucidée.

{\it La hauteur de Faltings (ou {hauteur modulaire})} $\rh_{\rm Fal}(\cdot)$ est définie de la 
manière suivante (cf. \cite[Par. 3]{Faltings-Endlichkeit}). Soit ${\omega}_{\CN(C)}$ la restriction 
de $\det(\Omega_{\CN(C)/\Spec\ \CO_{K_0}})$ à $\Spec\ \CO_{K_0}$ via 
la section nulle. Sur $\CN(C)_{K_0}=\Jac(C)$, on dispose 
d'une identification naturelle 
${\omega}_{\CN(C),K_0}\simeq\rH^0(\Jac(C),\Omega^g_{\Jac(C)/{K_0}}).$ 
Munissons 
${\omega}_{\CN(C)}$ de la métrique hermitienne donnée sur 
chaque composante $\CN(C)_\iota$ de $\CN(C)_\mC$ par la formule 
$$
\langle \alpha,\beta\rangle:={i^{g^2}\over(2\pi)^g}\int_{\CN(C)_\iota(\mC)}a\wedge\mtr{\beta}.
$$
Alors 
$$
\rh_{\rm Fal}(\Jac(C)_{\mtr{K}_0})=[K_0:\mQ]^{-1}\adeg(\mtr{\omega}_{\CN(C)}).
$$
Comme l'indique la notation, la quantité $\rh_{\rm Fal}(\Jac(C)_{\mtr{K}_0})$ 
ne dépend que de $C_{\mtr{K}_0}$ et non de $\CN(C)$. 

{\it La constante $\Theta_{\rm Max}(C_{\mtr{K}_0})$} est un élément de $\mR^*_+$ 
qui ne dépend que des points complexes de $C$ pour les divers 
plongements de $K_0$ dans $\mC$. 
Elle est définie de la manière suivante. 
Soit $M$ une surface 
de Riemann compacte de genre $g$ (donc l'ensemble des points complexes d'une courbe algébrique 
lisse et projective de genre $g$ sur $\mC$). Fixons 
$a_1,\dots a_g,b_{1},\dots,b_g$ une base symplectique de $\rH_1(M,\mZ)$. 
Soit $\omega_1,\dots,\omega_g$ la base 
de l'espace des $1$-formes  holomorphes $\rH^0(M,\omega_M)$ de $M$, 
qui est duale à $a_1,\dots,a_g$ pour l'accouplement 
$\int_{\gamma}\eta$ ($\gamma\in\rH_1(M,\mC)$, $\eta\in\rH^0(M,\omega_M)$). 
Soit $\tau$ la matrice $g\times g$ $$\tau:=\Big[\int_{b_j}\omega_i\Big]_{i,j}.$$
La fonction $\theta$ de Riemann (qui dépend du choix de la base 
symplectique) est alors par définition 
la fonction sur $\mC^g$
$$
\theta(z)=\theta(z,\tau):=\sum_{m\in\mZ^g}
\exp(2i\pi({1\over 2}\,^t m\tau m+\,^t m z))
$$
où l'on a identifié les élément de $\mZ^g$ à des matrices colonne. 
Soit $\Theta$ le lieu des zéros de $\theta$; c'est une sous variété analytique complexe fermée de $\mC^g$. Le lieu $\Theta$ est invariant sous 
l'action par translation des vecteurs colonne de la matrice  
$[\Id,\tau]$, qui a $g$ lignes et $2g$ colonnes. 
Le lieu $\Theta$ définit ainsi 
une sous-variété analytique complexe fermée du quotient $\mC^g/[\Id,\tau]$ 
de $\mC^g$ par le réseau engendré par les vecteurs colonnes de 
$[\Id,\tau]$. 
On le note  aussi $\Theta$ par abus de notation.  
Le quotient $\mC^g/[\Id,\tau]$ (qui est un tore complexe) 
s'identifie à la variété jacobienne $\Jac(M)$ de $M$. 
Il existe une unique métrique $\langle\cdot,\cdot\rangle$ sur $\CO(\Theta)$ 
telle que la forme de courbure associée à $(\CO(\Theta),\langle\cdot,\cdot\rangle)$ soit 
invariante par translation et 
telle que 
$$
\int_{\Jac(M)}\langle s,s\rangle\ \rd\mu=2^{-g/2},
$$
où $\rd\mu$ est la mesure de Haar sur $\Jac(M)$ de masse totale $1$ 
et $s$ est la  section canonique de $\CO(\Theta)$. 
Moret-Bailly (\cite[Par. 3.2, (3.2.2)]{Moret-Bailly-Sur}) démontre la formule explicite
$$
\langle s(z),s(z)\rangle:=\det(\Im(\tau)^{1\over 2}
\exp(-2\pi^t y(\Im(\tau))^{-1}y)|\theta(x+iy,\tau)|^2,
$$
où $z=x+iy\in\mC^g$. On identifie ici $z$ avec son image 
dans $\mC^g/[\Id,\tau]$. 
On définit enfin 
$$
\Theta_{\rm Max}(M):=\sup_{t\in\Jac(M)}\sqrt{\langle s(t),s(t)\rangle}.
$$
Cette constante n'est pas infinie, au vu de  la compacité de $\Jac(M)$.   
Enfin, on définit
$$
\Theta_{\rm Max}(C_{\mtr{K}}):=\sup_{\iota:K\hookrightarrow\mC}\Theta_{\rm Max}(C_\iota(\mC)).
$$
On constatera que la constante $\Theta_{\rm Max}(C_{\mtr{K}})$ ne 
dépend que de $C_{\mtr{K}}$. 

{\it Le cas des courbes de genre $2$}. 
On remarquera que si $g=2$ et que $P_0$ est fixé par une involution 
hyperelliptique de $C_{\mtr{K}_0}$, alors $\NT(j(\omega_C))=0$ car 
$j(\omega_C)$ est alors un point de $2$-torsion. Dans 
la même situation, la condition sur $(g-1)T$ est aussi superflue, 
car alors $$\Exc(W_{g-1})(\mtr{K}_0)
\cup[-1]^*\Exc(W_{g-1})(\mtr{K}_0)\subseteq j(C(\mtr{K}_0)).$$  

Le Théorème \ref{mainth} et  la Proposition \ref{mainprop}  
sont démontrées dans la section 3, la section 2 étant réservée a des préliminaires. 
On démontrera en fait dans la section 3 un théorème un peu plus fort 
que le Théorème \ref{mainth} (le Théorème 
\ref{mainthp}) mais dont l'énoncé est plus abscons. 
La section 3 traite l'exemple de la courbe $y^2+y=x^5$. Il y est fait 
un usage essentiel des calculs faits dans l'article \cite{Bost-Mestre-Calcul}.

{\bf Remerciements.} Je remercie V. Maillot pour nombre de 
discussions intéressantes et J. Boxall pour ses remarques et pour avoir 
signalé une erreur de calcul dans une version antérieure de cet article. Mes remerciements vont aussi à L. Moret-Bailly 
pour sa très utile monographie \cite{Moret-Bailly-Pinceaux}, sans laquelle 
cet article n'aurait pas vu le jour. Enfin, je suis reconnaissant à 
J.-F. Mestre, pour plusieurs calculs qu'il a fait dans le contexte du présent 
article et pour avoir partagé avec moi sa connaissance des fonctions  
thêta. 

\section{Préliminaires}

Dans cette section, on rappelle les cas particuliers des 
résultats de Moret-Bailly et Buium dont nous aurons besoin. 

\subsection{Fibrés en droites cubistes, selon L. Breen et L. Moret-Bailly}

Soit $L$ un corps de nombre et $B:=\Spec\ \CO_L$. 

Soit maintenant $\pi:\CA\to B$ un schéma en groupes commutatif 
lisse et de type fini sur $B$,  tel que $\CA_L$ est une variété abélienne. Soit 
$u:B\to\CA$ la section nulle. 

Soit $\mtr{\CL}$ un fibré en droites hermitien sur $\CA$. 
Nous dirons que $\mtr{\CL}$ est {\it cubiste}, si les conditions 
suivantes sont remplies. 
On requiert que $\CL$ soit cubiste (cf. \cite[chap. I]{Moret-Bailly-Pinceaux}
 ou \cite[chap. 2]{Breen-Fonctions}) pour cette notion), que $\c1(\mtr{\CL})_\mC)$ soit une forme différentielle invariante par translation et enfin que l'isomorphisme $u^*\CL\simeq\CO_B$ obtenu de 
la structure cubiste de $\CL$ soit une isométrie. On a muni ici le terme 
$\CO_B$ de la métrique triviale. On rappelle que l'on appelle 
{\it rigidification } (de $\CL$) un isomorphisme $u^*\CL\simeq\CO_B$. 

Pour chaque $P\in B$, soit $c_P$ l'ordre du groupe des composantes 
connexes de la fibres de $\CA$ au-dessus de $P$. Soit 
$c:={\rm ppcm}\{c_P|P\in B\}$. 
\begin{prop}[Moret-Bailly]
Soit $\CM$ un fibré cubiste sur $\CA_L$. 
Il existe alors une unique extension du fibré en droites 
cubiste $\CM^{\otimes 2c}$ à un fibré en droite cubiste sur $\CA$. 
\label{MB0}
\end{prop}
Pour la preuve, voir \cite[II, 1.2.1]{Moret-Bailly-Pinceaux}. 

\begin{prop}[Grothendieck; Breen; Moret-Bailly]
On suppose que $\CA$ est semi-stable sur $B$ et que 
les fibres de $\CA$ sur $B$ sont connexes. Alors 
le foncteur d'oubli de la catégorie des fibrés en droites 
cubistes sur $\CA$ vers la catégorie des fibrés en droites 
rigidifiés sur $\CA$ est une équivalence de catégorie.
\label{MB00}
\end{prop}
Dans la Proposition \ref{MB00}, les morphismes 
de la catégorie des fibrés en droites cubistes  
(resp. des fibrés en droites rigidifiés) sont les isomorphismes 
des fibrés en droites qui respectent les structures cubistes (resp. les
rigidifications). 
Pour la démonstration de la Proposition \ref{MB00}, voir par exemple \cite[I, Par. 2.6]{Moret-Bailly-Pinceaux}. 

\begin{theor}[Moret-Bailly]
Soit $\mtr{\CL}$ un fibré en droites hermitien cubiste sur $\CA$. 
Soit $z:B\to\CA$ une section. Alors
$[L:\mQ]^{-1}\adeg(z^*\mtr{\CL})$ est la hauteur de Néron-Tate 
associé au fibré en droite $\CL_L$ et au point 
$z_L\in\CA(L)$. 
\label{MB1}
\end{theor}
Pour la preuve, voir \cite[III, 4.4.1]{Moret-Bailly-Pinceaux}.

Enfin, nous allons citer le cas particulier de la  \og formule clé\fg\ 
de Moret-Bailly dont nous ferons usage.

On suppose à partir de maintenant que $\CA$ est un schéma abélien 
sur $B$. Soit $g$ la dimension relative de $\CA$ sur $B$. Soit $D\hookrightarrow\CA$ un diviseur effectif, symétrique 
et tel que $\CO(D)$ donne lieu à une polarisation principale sur 
chaque fibre de $\CA$ au-dessus de $B$. 

Nous  considérons $\CO_L$ 
comme un anneau arithmétique, que nous munissons comme tel de tous  ses plongements dans $\mC$. 

On munit $\CO(D)$  de l'unique 
métrique $\langle\cdot,\cdot\rangle$ dont la forme de courbure est invariante par translation et telle que pour tout $\iota\in\rHom(L,\mC)$, on ait
$$
\int_{\CA_\iota}\langle s,s\rangle\ \rd\mu_\iota=2^{-g/2},
$$
où $\rd\mu_\iota$ est la mesure de Haar sur $\CA_\iota(\mC)$ de masse totale $1$ 
et $s$ est la  section canonique de $\CO(D)$ (voir l'introduction). 
Soit $\mtr{\CO}(D)$ le fibré hermitien résultant.

\begin{theor}[Moret-Bailly]
L'égalité
$$
[L:\mQ]^{-1}\adeg(\mtr{\CO}(D)_0)={1\over 2}\rh_{\rm Fal}(\CA_{\mtr{L}})+{g\over 4}\log(4\pi)
$$
est vérifiée.
\label{MB2}
\end{theor}
Ici $\mtr{\CO}(D)_0$ est la restriction de $\mtr{\CO}(D)$ à 
$\Spec\ \CO_L$ via la section nulle.

Pour la preuve, voir \cite{Moret-Bailly-Sur}. 

\subsection{La géométrie modulo $l^2$ des courbes 
de genre $\geqslant 2$,  selon A. Buium}

Soit $R$ un anneau de valuation discrète, 
absolument non-ramifié et de corps résiduel 
de caractéristique $l$. Soit $F$ son corps de fractions. 

Soit $\wt{E}$ une courbe propre et lisse sur 
$\Spec\ R$. On suppose que les fibres géométriques 
de $\wt{E}$ sont de genre $g\geqslant 2$. On suppose 
donné un point $Q\in\wt{E}(F)$. La propriété universelle des 
modèles de Néron implique que le plongement 
jacobien $\wt{E}_F\hookrightarrow\Jac(\wt{E}_F)$ s'étend 
en un plongement $\wt{E}\to\CN(\wt{E}_F)$ dans le modèle de 
Néron $\CN(\wt{E}_F)$ de $\Jac(\wt{E}_F)$ sur $R$. 

Le théorème ci-dessous a été établi par Buium dans 
le cadre de sa  démonstration effective de la conjecture de Manin-Mumford pour les courbes: 

\begin{theor}[Buium]
L'inégalité
$$
\#\wt{E}(R/l^2 R)\cap l\cdot\CN(\wt{E}_F)(R/l^2 R)\leqslant\Bu_l
$$
est vérifiée.
\label{BU1}
\end{theor}
Pour la démonstration, voir \cite[fin de la preuve du Th. 1.11]{Buium-p-jets}.

\section{Démonstration de \ref{mainth} et \ref{mainprop}}

Nous allons d'abord démontrer le Théorème \ref{mainth}. 

Ce théorème est impliqué par le Théorème \ref{mainthp} ci-dessous. 

Soit $K$ une extension finie de $K_0$  et 
$\mfp$ un idéal premier de $\CO_K$.  Soit $\mfp_0:=\mfp\cap\CO_{K_0}$ 
et soit $p$ le générateur positif 
de $\mfp\cap\mZ$. Soit  $q:=\#\CO_{K_0}/\mfp_0$. 

\begin{theor}
Soit $T\in\Tor(\Jac({C})(K))$. On suppose que 
\begin{description} 
\item[\ (1)] la courbe $C$ a réduction semi-stable sur $K_0$;
\item[\ (2)] $p>2$; 
\item[\ (3)] $\mfp_0$ est une place de bonne réduction de $C$;
\item[\ (4)] $\mfp_0$ est non-ramifié sur $\mZ$;
\item[\ (5)] $(\ordre(T),p)=1$;
\item[\ (6)] la section $\Spec\ \CO_{K}\to\CN(C_K)$ correspondant 
à $T$ se factorise par la composante neutre de $\CN(C_K)$.
\end{description}
Alors l'inégalité
\begin{eqnarray}
d_\mfp(T,{C})&\geqslant&p^{-\textstyle\Big(1+ 2\cdot\rL_q\cdot[K_0:\mQ]\cdot\Big|\log\Theta_{\rm Max}(C_{\mtr{K}})+[K_0:\mQ]^{-1}\adeg(\mtr{\CO}(\Zar(W_{g-1}))_0)\Big|\Big)}\ \ \ \ \ \ 
\label{mform}
\end{eqnarray}
est vérifiée chaque fois que 
$$
(g-1)T\not\in\Exc(W_{g-1})(\mtr{K})\cup[-1]^*\Exc(W_{g-1})(\mtr{K}).
$$
\label{mainthp}
\end{theor}
\beginProof 
Puisque $\CN(C_K)_{\mfp}$ est un schéma en groupes sur $\CO_{\mfp}$, on a 
$$
\rv_\mfp(T,{C}_{K})\leqslant \rv_\mfp((g-1)T,W_{g-1,K})
$$
et on est donc ramené à estimer $\rv_\mfp((g-1)T,W_{g-1,K})$. 

Soit $W':=W_{g-1,K}\cup[-1]^*W_{g-1,K}$. 
Soit $T':=(g-1)T$. 
Soit 
$t:\Spec\ \CO_K\to\CN(C_K)$ la section correspondante à 
$T'$. 

Si $\mfl$ est un idéal premier de 
$\CO_{K}$, on écrira $\rN_\mfl:=\#\CO_{K}/(\mfl\cap\mZ)$. 

Soit $\mtr{L}:=\mtr{\CO}(\Zar(W_{g-1}))$. Soit 
 $\mtr{L}_0$ pour la restriction de $\mtr{L}$ à 
$\Spec\ \CO_{K_0}$ via la section nulle. 
Soit $s_L$ la section canonique de $L$ s'annulant sur $\Zar(W_{g-1})$ 
et $\langle\cdot,\cdot\rangle_L$ la métrique de $\mtr{L}$. 

Soit $$\mtr{M}:=\mtr{L}\otimes[-1]^*\mtr{L}\otimes\mtr{L}^{\vee,\otimes 2}_0.$$ 

On peut supposer sans restriction de généralité que 
$K$ est la clôture galoisienne du corps de définition 
de $T$ sur $K_0$. 

On calcule
\begin{eqnarray*}
v_{\mfp}(T',{W_{g-1,K}})&\leqslant& v_{\mfp}(T',W')\stackrel{\rm (a)}{\leqslant} 
\sum_{\mfl\in\Spec\ \CO_{K}}\log(\rN_\mfl)\cdot v_\mfl(T',W')\\
&\stackrel{\rm (b)}{=}&(2\mfn_K)^{-1}\adeg\big(t^*(\mtr{L}\otimes 
[-1]^*\mtr{L})^{\otimes 2\mfn_K}\big)+{1\over 2}\sum_{\iota:K\hookrightarrow\mC}\log
\langle s_L(T'),s_L(T')\rangle_L\\
&&+ {1\over 2}\sum_{\iota:K\hookrightarrow\mC}\log
\langle s_L(-T'),s_L(-T')\rangle_L\\
&\stackrel{\rm (c)}{\leqslant}&(2\mfn_K)^{-1}\adeg(t^*\mtr{M}^{\otimes 2\mfn_K})+{2[K:\mQ]}\log\Theta_{\rm Max}(C_{\mtr{K}})+
2[K:K_0]\adeg(\mtr{L}_0)\\
&\stackrel{\rm (d)}{=}&{2[K:\mQ]}\log\Theta_{\rm Max}(C_{\mtr{K}})+
2[K:K_0]\adeg(\mtr{L}_0).\\
\end{eqnarray*}

L'inégalité (a) est justifiée par la condition (2) du Théorème \ref{mainthp}. 
L'égalité (b) est justifié par la définition du degré arithmétique $\adeg(\cdot)$. 
L'égalité (c) est justifiée par le fait que
$$
[K:K_0]\adeg(\mtr{L}_0)=\adeg(\mtr{L}_{\CO_K,0}).
$$
En effet, les composantes neutres de $\CN(C)_K$ et de 
$\CN(C_K)$ coïncident à cause de la condition de semi-stabilité (1) 
du Théorème \ref{mainthp}; par ailleurs, la clôture de Zariski 
commute aux changements de base plat (cf. \cite[IV, 2.4.11]{EGA}). 
L'égalité (d) est justifiée les Théorèmes \ref{MB0}, \ref{MB00} et \ref{MB1} 
ainsi que la condition (6) du Théorème \ref{mainthp}. 

Nous allons maintenant estimer la quantité $[K:\mQ]$ au moyen des résultats de Buium. 
On remarque tout d'abord que la condition (6) du Théorème \ref{mainthp} est 
invariante par restriction ou extension du corps $K$. Ceci est à nouveau une conséquence de la condition de semi-stabilité, qui assure que les composantes neutres de $\CN(C)_K$ et de 
$\CN(C_K)$ coïncident. On peut donc sans restriction de la généralité supposer  
que $K$ est le corps engendré au-dessus de $K_0$ par les corps 
de définition des points de $N$-torsion 
de $\Jac(C)(\mtr{K}_0)$, où $N:=\ordre(T)$. Ce corps est par construction 
galoisien sur $K_0$; par ailleurs, comme $(N,p)=1$, il est non-ramifié 
au-dessus de $\mfp_0$.  On peut aussi supposer sans restriction de la généralité 
que $\rv_\mfp(T,C_K)>1$, vu la forme de la partie droite 
de l'inégalité \refeq{mform} du Théorème \ref{mainthp}. 

Soit
$$
E:=C(\CO_K/\mfp^2)\cap p\cdot\CN(C)(\CO_K/\mfp^2).
$$
Soit $\mF:=\CO_K/\mfp$. Le corps  $\mF$ est par construction une extension 
de $\CO_{K_0}/\mfp_0$ et l'on identifie ce dernier corps 
à $\mF_q$. Soit $G_\mfp\subseteq\Gal(K|K_0)$ le groupe 
de décomposition de $\mfp$. Les conditions (5) et (3) du Théorème \ref{mainthp} 
assurent que $K$ est non-ramifié au-dessus de $\mfp_0$ et on dispose 
donc d'une identification canonique $G_\mfp\simeq\Gal(\mF|\mF_q)$. 
L'ensemble $E$ est naturellement $G_\mfp$ invariant et l'application 
naturelle $\phi:E\to\CN(C)(\mF)$ est $G_\mfp$-équivariante. 
On remarque maintenant que le Théorème \ref{BU1} implique 
l'inégalité $\#E\leqslant\Bu_p$ (on a utilisé ici les conditions 
(3) et (4) du Théorème \ref{mainthp}). On en déduit que pour tout 
$e\in E$, on a 
$$
\deg_{\mF_q}(\phi(e))\leqslant\Bu_p.
$$
On remarque qu'il y 
des flèches naturelles injectives
$$
\Tor^p(\CN(C)(K))=\Tor^p(\CN(C)(\CO_K))\hookrightarrow\Tor^p(\CN(C)(\CO_K/\mfp^2))\hookrightarrow 
\Tor^p(\CN(C)(\mF)).
$$
Puisque $\Tor^p(\CN(C)(K))\subseteq p\cdot\CN(C)(K)$ 
et $T\ (\mod\ \mfp^2)\in C(\CO_K/\mfp^2)$, 
on voit que 
$$
\deg_{\mF_q}(T\ (\mod\ \mfp))=:d\leqslant\Bu_p.
$$ 
Les estimées de Hasse-Weil (cf. par ex. \cite[Par. 19, Th. 19.1]{Milne-Abelian}) impliquent que 
$$
\#\CN(C)(\mF_{q^{d}})\leqslant q^{d\cdot g}+(2^{2g}-2g-1)q^{d(g-1)}+2gq^{d(g-1/2)}.
$$
Rappelons que $N=\ordre(T)$. 
Puisque $N|\#\CN(C)(\mF_{q^{d}})$, on déduit que 
$$
N\leqslant q^{\Bu_p\cdot g}+(2^{2g}-2g-1)q^{\Bu_p(g-1)}+2gq^{\Bu_p(g-1/2)}.
$$
Remarquons maintenant que par construction 
$$
[K:K_0]\leqslant\#\GL_{2g}(\mZ/N\mZ)\leqslant N^{4g^2}.
$$
On obtient pour finir 
$$
[K:K_0]\leqslant \big[q^{\Bu_p\cdot g}+(2^{2g}-2g-1)q^{\Bu_p(g-1)}+2gq^{\Bu_p(g-1/2)}\big]^{4g^2}=:L_q.\qed
$$
{\it Ceci termine la démonstration du Théorème \ref{mainth}.}

Nous allons maintenant passer à la démonstration de la Proposition \ref{mainprop}.

On suppose que $C_{K_0}$ a bonne réduction partout. 

On rappelle qu'il existe un point $\kappa\in\Jac(C)(\mtr{K}_0)$ 
telle que $W_{g-1,\mtr{K}_0}+\kappa$ est un diviseur $\Theta$. En particulier, 
$W_{g-1,\mtr{K}_0}+\kappa$ est alors symétrique. 
Comme l'égalité de la Proposition \ref{mainprop} est invariante par 
extension du corps $K_0$, on peut supposer sans restriction de généralité 
que $\kappa$ est défini sur $K_0$. On peut montrer que 
$-2\kappa=j(\omega_C)$ (cf. \cite[I, 5.]{Arbarello-Geometry}). Soit $k:\Spec\ \CO_{K_0}\to\CN(C)$ 
la section correspondant à $\kappa$. 
On peut maintenant calculer
\begin{eqnarray*}
[K_0:\mQ]^{-1}\adeg(\mtr{\CO}(\Zar(W_{g-1}))_0)&=&
[K_0:\mQ]^{-1}\adeg(k^*\mtr{\CO}(\Zar(W_{g-1}+\kappa)))\\
&=&
[K_0:\mQ]^{-1}\adeg(k^*\mtr{\CO}(\Zar(W_{g-1}+\kappa))\otimes
\mtr{\CO}(\Zar(W_{g-1}+\kappa))_0^\vee)\\
&&+[K_0:\mQ]^{-1}\adeg(\mtr{\CO}(\Zar(W_{g-1}+\kappa))_0)\\
&\stackrel{\rm (a)}{=}&
\NT(\kappa)+{1\over 2}\rh_{\rm Fal}(\Jac(C)_{\mtr{K}_0})+{g\over 4}\log(4\pi)\\
&\stackrel{\rm (b)}{=}&
{1\over 4}\NT(j(\omega_C))+{1\over 2}\rh_{\rm Fal}(\Jac(C)_{\mtr{K}_0})+{g\over 4}\log(4\pi).
\end{eqnarray*}
L'égalité (a) est justifiée par le Théorème \ref{MB1} et par 
le Théorème \ref{MB2}. L'égalité (b) est justifiée par le fait que 
la hauteur de Néron-Tate est quadratique. 

{\it Ceci  conclut la démonstration de la Proposition \ref{mainprop}.}

\section{La courbe $y^2+y=x^5$}

Nous nous servons dans cette section 
des calculs fait dans \cite{Bost-Mestre-Calcul} pour obtenir 
une borne explicite pour la distance $p$-adique dans le cas 
de la courbe $C$ donné en coordonnée affine sur $\mQ$ par 
l'équation $y^2+y=x^5$. Soit 
$$
K_0:=\mQ(\sqrt{1-\exp(2i\pi/5)},\sqrt[5]{2}).
$$
On vérifie que $[K_0:\mQ]=40$ et que seuls $2$ et $5$ se ramifient dans $K_0$. 
Il est montré dans \cite[Par. 4.1]{Bost-Mestre-Calcul} que $C$ est de genre 
$2$ et qu'elle acquiert bonne réduction partout sur $K_0$. 
Soit $\zeta_5:=\exp(2i\pi/5)$. 
On 
montre dans \cite[Par. 4.4]{Bost-Mestre-Calcul} qu'il existe une base symplectique de $H_1(C(\mC),\mZ)$ telle que 
\begin{displaymath}
\tau=
\left(\begin{array}{cc}
-\zeta^4_5 & \zeta^2_5+1\\
\zeta^2_5+1 & \zeta^2_5-\zeta^3_5
\end{array}
\right).
\end{displaymath}
Donc
$$
\Theta_{\rm Max}=
\sup_{x+iy\in\mC^2}
\det(\Im(\tau))^{1\over 4}
\exp(-\pi^t y(\Im(\tau))^{-1}y)\big|\sum_{m\in\mZ^2}
\exp(2i\pi({1\over 2}\,^t m\tau m+\,^t m z))\big|.
$$
La constante $\Theta_{\rm Max}$ peut être approximée par ordinateur. 
J.-F. Mestre a eu la gentillesse de me communiquer l'approximation 
suivante, réalisée avec le programme MAPLE:
$$
\Theta_{\rm Max}\approx 1.06639277369136206671054075
$$
Par ailleurs, on a
$$
\rh_{\rm Fal}(C_{\mtr{\mQ}})=2\log(2\pi)-
{1\over 2}\log\big(\Gamma(1/5)^5\Gamma(2/5)^3\Gamma(3/5)\Gamma(4/5)^{-1}\big)
\approx
-1.452509239645644650317707042
$$
Pour cela, voir \cite[Par. 4.4, Prop. 12]{Bost-Mestre-Calcul}. 
Ainsi
$$
\log\Theta_{\rm Max}(C_{\mtr{K}})+{1\over 2}\rh_{\rm Fal}(C_{\mtr{\mQ}}))
+{1\over 2}
\log(4\pi)
\approx 0.60353921716278764047528474
$$
J.-F. Mestre me communique également que le nombre ci-dessus 
ressemble à s'y méprendre au nombre
$
{3\over 8}\log(5)$. 
Après le changement de variable $z=2y+1$ et $t=\sqrt[5]{4}\cdot x$, la courbe 
$C_{K_0}$ admet l'équation affine
$$
z^2=t^5+1.
$$
Soit $p$ non-ramifié dans $K_0$ (i.e. $p\not\in\{2,5\}$). 
Soit $P_0$ le point donné par les coordonnées $z=1$, $t=0$. Le point 
$P_0$ est invariant par l'involution hyperelliptique $z\to -z$, $t\to t$.

Soit $C_p$ la courbe 
obtenue sur $\mtr{\mQ}_p$ à partir de $C$. Soit $C_p\hookrightarrow 
\Jac(C_p)$ le plongement jacobien obtenu à partir de $j_{P_0}$. 

 Le Théorème \ref{mainthp} et la Proposition \ref{mainprop} impliquent maintenant
que 
\begin{eqnarray*}
d_p(T,{C}_p)&\geqslant&p^{\textstyle\Big(1+ 80\cdot\rL_{q}\cdot\Big|
\log\Theta_{\rm Max}(C_{\mtr{K}})+{1\over 2}\rh_{\rm Fal}(C_{\mtr{\mQ}}))+{1\over 2}
\log(4\pi)\Big|\Big)}
\end{eqnarray*}
pour tout point de torsion $T\in\Jac(C_p)(\mtr{\mQ}_p)$ telle que 
$(\ordre(T),p)=1$ et tel que $T\not\in C_p(\mtr{\mQ}_p)$. 

On rappelle que $q:=\#\CO_{K_0}/\mfp_0$. On a ainsi 
l'estimée $q\leqslant p^{40}$. On remarque qu'il est préférable 
de choisir des nombres premiers $p$ tels que 
$q$ est petit, car la fonction $\rL_q$ est doublement exponentielle en $p$.

\end{document}